\documentclass{amsart}
\usepackage{setspace}
\usepackage{a4}
\usepackage{amsthm}
\usepackage{latexsym}
\usepackage{amsfonts}
\usepackage{graphicx}
\usepackage{textcomp}
\usepackage{cite}
\usepackage{enumerate}
\usepackage{amssymb}
\usepackage{hyperref}
\usepackage{amsmath}
\usepackage{tikz}
\usepackage[mathscr]{euscript}
\usepackage{mathtools}
\newtheorem{theorem}{Theorem}[section]

\newtheorem{definition}[theorem]{Definition}
\newtheorem{example}[theorem]{Example}

\newtheorem{problem}[theorem]{Problem}

\title{This is the title}
\usepackage{fancyhdr}

\begin{document}
\hrule\hrule\hrule\hrule\hrule
\vspace{0.3cm}	
\begin{center}
{\bf{p-adic Grothendieck Inequality, p-adic Johnson-Lindenstrauss Flattening and p-adic Bourgain-Tzafriri Restricted Invertibility Problems}}\\
\vspace{0.3cm}
\hrule\hrule\hrule\hrule\hrule
\vspace{0.3cm}
\textbf{K. Mahesh Krishna}\\
School of Mathematics and Natural Sciences\\
Chanakya University Global Campus\\
NH-648, Haraluru Village\\
Devanahalli Taluk, 	Bengaluru  North District\\
Karnataka  562 110 India\\
Email: kmaheshak@gmail.com\\

Date: \today
\end{center}

\hrule\hrule
\vspace{0.5cm}
%--------------------------------------
\textbf{Abstract}: We formulate p-adic versions of following three: (1)  Grothendieck Inequality, (2) Johnson-Lindenstrauss Flattening Lemma, (3) Bourgain-Tzafriri Restricted Invertibility Theorem.

\textbf{Keywords}:   Grothendieck Inequality, Johnson-Lindenstrauss  Lemma, Bourgain-Tzafriri Restricted Invertibility, p-adic Hilbert space.

\textbf{Mathematics Subject Classification (2020)}: 46S10.\\

\hrule

\vspace{0.5cm}
%\tableofcontents

\section{p-adic Grothendieck Inequality Problem}
 In 1953, Grothendieck derived the following inequality which is one of the greatest inequalities of all time. 
  \begin{theorem}\cite{BLEI, FRIEDLANDLIMZHANG, RIETZ, GROTHENDIECK, ALBIACKALTON, DIESTELFOURIESWART, PISIERGROTHENDIECK, LINDENSTRAUSS}\label{G}
  	 \textbf{(Grothendieck Inequality)} There is a universal constant $K_G$ satisfying the following: For every Hilbert space $\mathcal{H}$ and for all $m, n \in \mathbb{N}$, if a scalar matrix $[a_{j,k}]_{1\leq j \leq m, 1\leq k \leq n}$ satisfies
  	\begin{align*}
  		\left|\sum_{j=1}^m \sum_{k=1}^na_{j,k}s_jt_k\right|\leq \left(\max_{1\leq j \leq m}|s_j|\right) \left(\max_{1\leq k \leq n}|t_k|\right),\quad  \forall s_j , t_k \in \mathbb{C}, 
  	\end{align*}
  then 
  \begin{align*}
  \left| \sum_{j=1}^m \sum_{k=1}^na_{j,k}\langle u_j, v_k\rangle\right|	\leq K_G\left(\max_{1\leq j \leq m}\|u_j\|\right)  \left(\max_{1\leq k \leq n}\|v_k\|\right), \quad  \forall u_j , v_k \in \mathcal{H}. 
  \end{align*}
  	 \end{theorem}
   We formulate a problem which is p-adic analogue of Theorem \ref{G}. We first recall the definition of p-adic Hilbert space.
   \begin{definition}\cite{DIAGANABOOK, DIAGANARAMAROSON, KALISCH} \label{PADICDEF}
   	Let $\mathbb{K}$ be a non-Archimedean  valued field (with valuation $|\cdot|$) and $\mathcal{X}$ be a non-Archimedean Banach space (with norm $\|\cdot\|$) over $\mathbb{K}$. We say that $\mathcal{X}$ is a \textbf{p-adic Hilbert space} if there is a map (called as p-adic inner product) $\langle \cdot, \cdot \rangle: \mathcal{X} \times \mathcal{X} \to \mathbb{K}$ satisfying following.
   	\begin{enumerate}[\upshape (i)]
   		\item If $x \in \mathcal{X}$ is such that $\langle x,y \rangle =0$ for all $y \in \mathcal{X}$, then $x=0$.
   		\item $\langle x, y \rangle =\langle y, x \rangle$ for all $x,y \in \mathcal{X}$.
   		\item $\langle x, \alpha y+z \rangle =\alpha \langle x,  y \rangle+\langle x,z\rangle$ for all  $\alpha  \in \mathbb{K}$, for all $x,y,z \in \mathcal{X}$.
   		\item $|\langle x, y \rangle |\leq \|x\|\|y\|$ for all $x,y \in \mathcal{X}$.
   	\end{enumerate}
   \end{definition}
Following is the  standard finite dimensional example.
\begin{example}\cite{KALISCH}
	Let $p$ be a prime. For $d \in \mathbb{N}$, let $\mathbb{Q}_p^d$ be the standard p-adic Hilbert space equipped with the inner product 
	\begin{align*}
		\langle (a_j)_{j=1}^d,(b_j)_{j=1}^d\rangle \coloneqq \sum_{j=1}^da_jb_j,  \quad \forall (a_j)_{j=1}^d,(b_j)_{j=1}^d \in \mathbb{Q}_p^d
	\end{align*}
	and norm 
	\begin{align*}
		\|(x_j)_{j=1}^d\|\coloneqq \max_{1\leq j \leq d}|x_j|, \quad \forall (x_j)_{j=1}^d\in 	\mathbb{Q}_p^d.
	\end{align*}
\end{example}
Based on Theorem \ref{G}, we formulate following problem.
   \begin{problem}
 	\textbf{(p-adic  Grothendieck Inequality Problem)} 	Let $\mathbb{K}$ be a  non-Archimedean  valued field. Whether there is a  universal constant $K_\mathbb{K}$  satisfying the following property: for every   p-adic Hilbert space $\mathcal{X}$  over  $\mathbb{K}$, for all $m,n \in \mathbb{N}$,	if $[a_{j,k}]_{1\leq j \leq m, 1\leq k \leq n}\in \mathbb{M}_{m\times n}(\mathbb{K})$ satisfies 
 		\begin{align*}
 		\left|\sum_{j=1}^m \sum_{k=1}^na_{j,k}s_jt_k\right|\leq \left(\max_{1\leq j \leq m}|s_j|\right) \left(\max_{1\leq k \leq n}|t_k|\right),\quad  \forall s_j , t_k \in \mathbb{K}, 
 	\end{align*}
 	then 
 	\begin{align*}
 		\left| \sum_{j=1}^m \sum_{k=1}^na_{j,k}\langle u_j, v_k\rangle\right| 	\leq K_\mathbb{K}\left(\max_{1\leq j \leq m}\|u_j\|\right)  \left(\max_{1\leq k \leq n}\|v_k\|\right), \quad  \forall u_j , v_k \in \mathcal{X}.
 	\end{align*}
 \end{problem}
 
 \begin{problem}
	\textbf{(p-adic  Grothendieck Inequality Problem)} Let $\mathbb{K}$ be a  non-Archimedean  valued field. Whether there is a universal constant $K_\mathbb{K}$  satisfying the following property:	for every   p-adic Hilbert space $\mathcal{X}$  over  $\mathbb{K}$, for all $m,n \in \mathbb{N}$, for every matrix  $[a_{j,k}]_{1\leq j \leq m, 1\leq k \leq n}\in \mathbb{M}_{m\times n}(\mathbb{K})$,  we have 
		\begin{align*}
		&\sup\left\{\left| \sum_{j=1}^m \sum_{k=1}^na_{j,k}\left \langle u_j, v_k\right \rangle\right| : u_1, \dots, u_m, v_1, \dots, v_n \in \mathcal{X}, \|u_1\|=\cdots=\|u_m\|=\|v_1\|=\cdots=\|v_n\|=1\right\} \leq \\
		&	K_\mathbb{K}\sup\left\{\left| \sum_{j=1}^m \sum_{k=1}^n a_{j,k}s_j t_k\right| : s_1, \dots, s_m, t_1, \dots, t_n \in \mathbb{K}, |s_1|=\cdots=|s_m|=|t_1|=\cdots =|t_n| =1\right\}.
		\end{align*}
\end{problem}

    \section{p-adic Johnson-Lindenstrauss Flattening Problem}
   Johnson and  Lindenstrauss showed that points in high dimensional Euclidean space can be mapped to low dimensional Euclidean space with little distortion to distances between every pair of points,  using Lipschitz maps
    \cite{JOHNSONLINDENSTRAUSS}.
\begin{theorem}\cite{MATOUSEKBOOK, JOHNSONLINDENSTRAUSS}\label{JL}
	\textbf{(Johnson-Lindenstrauss Flattening Lemma)} Let $M, N\in \mathbb{N}$ and  $\mathbf{x}_1, \mathbf{x}_2, \dots, \mathbf{x}_M\in \mathbb{R}^N$. For each $0<\varepsilon<1$, there exists a Lipschitz map $f: \mathbb{R}^N\to \mathbb{R}^m$  and a real $r>0$ such that 
		\begin{align*}
			r(1-\varepsilon)\|\mathbf{x}_j-\mathbf{x}_k\|\leq \|f(\mathbf{x}_j)-f(\mathbf{x}_k)\|\leq r(1+\varepsilon)\|\mathbf{x}_j-\mathbf{x}_k\|, \quad \forall 1 \leq j, k\leq M, 
		\end{align*}
		where 
		\begin{align*}
			m=O\left(\frac{\log M}{\varepsilon^2}\right).
	\end{align*}
\end{theorem}
Over the period, some improvements of Theorem \ref{JL} were obtained. We recall them. 
\begin{theorem}\cite{FRANKLMAEHARA, JOHNSONLINDENSTRAUSS} \textbf{(Johnson-Lindenstrauss  Flattening  Lemma: Frankl-Maehara form)}
		Let $0<\varepsilon<\frac{1}{2}$   and $ M\in \mathbb{N}$. Define 
		\begin{align*}
			m (\varepsilon, M)\coloneqq \left\lceil 9\frac{1}{\varepsilon^2-\frac{2\varepsilon^3}{3}}\log M \right \rceil +1.
		\end{align*} 
		If  $M> m (\varepsilon, M)$, then for any  $\mathbf{x}_1, \mathbf{x}_2, \dots, \mathbf{x}_M\in \mathbb{R}^M$,  there exists 	a map $f: \{\mathbf{x}_j\}_{j=1}^M\to \mathbb{R}^m$   such that 
		\begin{align*}
			(1-\varepsilon)\|\mathbf{x}_j-\mathbf{x}_k\|^2\leq \|f(\mathbf{x}_j)-f(\mathbf{x}_k)\|^2\leq (1+\varepsilon)\|\mathbf{x}_j-\mathbf{x}_k\|^2, \quad \forall 1 \leq j, k\leq M.
	\end{align*}
\end{theorem}
\begin{theorem}\cite{DASGUPTAGUPTA, JOHNSONLINDENSTRAUSS}
	\textbf{(Johnson-Lindenstrauss Flattening Lemma: Dasgupta-Gupta form)}
		Let $M, N\in \mathbb{N}$ and  $\mathbf{x}_1, \mathbf{x}_2, \dots, \mathbf{x}_M\in \mathbb{R}^N$. Let $0<\varepsilon<1$. Choose any natural number $m$ such that 
		\begin{align*}
			m>4\frac{1}{\frac{\varepsilon^2}{2}-\frac{\varepsilon^3}{3}}\log M. 
		\end{align*}
		Then there exists a  map $f: \mathbb{R}^N\to \mathbb{R}^m$   such that 
		\begin{align*}
			(1-\varepsilon)\|\mathbf{x}_j-\mathbf{x}_k\|^2\leq \|f(\mathbf{x}_j)-f(\mathbf{x}_k)\|^2\leq (1+\varepsilon)\|\mathbf{x}_j-\mathbf{x}_k\|^2, \quad \forall 1 \leq j, k\leq M.
		\end{align*}
		The map $f$ can be found in randomized polynomial time.
\end{theorem}
\begin{theorem} \cite{JOHNSONLINDENSTRAUSS, FOUCARTRAUHUT} \textbf{(Johnson-Lindenstrauss Flattening Lemma: matrix form)} \label{MF}
		There is a universal constant $C>0$ satisfying the following. Let $0<\varepsilon<1$, $M,N \in \mathbb{N}$ and $\mathbf{x}_1, \mathbf{x}_2, \dots, \mathbf{x}_M\in \mathbb{R}^N$. For each natural number 
		\begin{align*}
			m>\frac{C}{\varepsilon^2}\log M,
		\end{align*}
		there exists a matrix $M\in \mathbb{M}_{m\times N}(\mathbb{R})$ such that 
		\begin{align*}
			(1-\varepsilon)\|\mathbf{x}_j-\mathbf{x}_k\|^2\leq \|M(\mathbf{x}_j-\mathbf{x}_k)\|^2\leq (1+\varepsilon)\|\mathbf{x}_j-\mathbf{x}_k\|^2, \quad \forall 1 \leq j, k\leq M.
	\end{align*}
\end{theorem}
Based on Theorem \ref{MF}, we formulate following problem.
\begin{problem}\textbf{(p-adic Johnson-Lindenstrauss Flattening Problem)}
	Let $\mathbb{K}$ be a non-Archimedean valued field. What is the best function $\phi: (0,1)\times \mathbb{N} \to (0,\infty)$ satisfying the following: 	There is a universal constant $C>0$ (which may depend upon $\mathbb{K}$) satisfying the following. Let $0<\varepsilon<1$, $M,N \in \mathbb{N}$ and $\mathbf{x}_1, \mathbf{x}_2, \dots, \mathbf{x}_M\in \mathbb{K}^N$. For each natural number 
	\begin{align*}
		m>C\phi (\varepsilon, M),
	\end{align*}
	there exists a matrix $M\in \mathbb{M}_{m\times N}(\mathbb{K})$ such that 
	\begin{align*}
		(1-\varepsilon)\|\mathbf{x}_j-\mathbf{x}_k\| \leq  \|M(\mathbf{x}_j-\mathbf{x}_k)\| \leq (1+\varepsilon)\|\mathbf{x}_j-\mathbf{x}_k\|, \quad \forall 1 \leq j, k\leq M.
\end{align*}
\end{problem}

  \section{p-adic Bourgain-Tzafriri Restricted Invertibility Problem}
   Let $\mathbb{R}^d$ be the standard $d$-dimensional Euclidean  space with canonical orthonormal basis $\{e_j\}_{j=1}^d$.  Bourgain and  Tzafriri  in 1987 showed that any matrix can be inverted restrictively with control over norm   \cite{BOURGAINTZAFRIRI}.
  \begin{theorem}\cite{BOURGAINTZAFRIRI}\label{BOURGAINTZARIRITHEOREM}  \textbf{(Bourgain-Tzafriri Restricted Invertibility Theorem)}
  		There are universal constants $A>0$, $c>0$ satisfying the following property. If  $d \in \mathbb{N}$, and  $T:\mathbb{R}^d \to \mathbb{R}^d$ is a linear operator   with  $\|Te_j\|=1$, $\forall 1\leq j\leq d$, then there exists a subset $\sigma \subseteq \{1, \dots, d\}$ of cardinality 
  		\begin{align*}
  			\operatorname{Card}(\sigma)\geq \frac{cd}{\|T\|^2}
  		\end{align*}
  		such that 
  		\begin{align*}
  			\left\|\sum_{j \in \sigma}a_jTe_j\right\|^2\geq A \sum_{j \in \sigma}|a_j|^2, \quad \forall a_j \in  \mathbb{R}, \forall j \in \sigma.
  	\end{align*}
  \end{theorem}  
  Recently,   Spielman and Srivastava gave a simple proof of Theorem \ref{BOURGAINTZARIRITHEOREM} by proving a generalization of  Theorem \ref{BOURGAINTZARIRITHEOREM}  \cite{SPIELMANSRIVASTAVA} due to Vershynin \cite{VERSHYNIN}.
Based on Theorem \ref{BOURGAINTZARIRITHEOREM}, we formulate following problem.
  \begin{problem}
  	\textbf{(p-adic Bourgain-Tzafriri Restricted Invertibility Problem)}
  	Let $\mathbb{K}$ be a non-Archimedean valued field. 		Whether there are universal constants $A>0$, $c>0$ (which may depend upon $\mathbb{K}$) satisfying the following property: If  $d \in \mathbb{N}$, and  $T:\mathbb{K}^d \to \mathbb{K}^d$ is a linear operator   with  $\|Te_j\|=1$, $\forall 1\leq j\leq d$, then there exists a subset $\sigma \subseteq \{1, \dots, d\}$ of cardinality 
  		\begin{align*}
  			\operatorname{Card}(\sigma)\geq \frac{cd}{\|T\|^2}
  		\end{align*}
  		such that 
  		\begin{align*}
  			\left\|\sum_{j \in \sigma}a_jTe_j\right\|^2\geq A \max_{j \in \sigma}|a_j|^2, \quad \forall a_j \in  \mathbb{K}, \forall j \in \sigma.
  	\end{align*}
  \end{problem}

 \bibliographystyle{plain}
 \bibliography{reference.bib}

\end{document}